# On a congruence only holding for primes II

E. Vantieghem[*]

**Abstract**

We present a primality criterium based on congruences for cyclotomic polynomials, and point out a way to generalize our result in order to obtain a family of similar criteria. No practical use is aimed however.

In this note we propose a family of congruences each of which holds if and only if a natural number $m$ greater than $2$ is prime. The flavour of our result is the same as in [1] and [2], so that we cannot expect any practical use. First we treat the simplest case, because it requires a different argument than the others.

THEOREM 1. Let $m$ be a natural number greater than 2. Then $m$ is prime if and only if
$$\prod_{j=1}^{m-1}(2^j+1) \equiv 1 \pmod{2^m-1}. \qquad (1)$$

Our proof is based on some facts from the theory of cyclotomic polynomials.

We recall that the $n$-th cyclotomic polynomial $\Phi_n(X)$ is defined by
$$\Phi_n(X) = \prod_{\substack{(k,n)=1 \\ 1 \le k < n}} (X - \zeta^k),$$

where $\zeta$ is any primitive $n$-th root of unity. Further we shall use the identity :
$$X^n - 1 = \prod_{d \mid n} \Phi_d(X) \qquad (*)$$

and the fact that the coefficients of $\Phi_n(X)$ are integers.

We will also need the following congruence in $\mathbb{Z}[X,Y]$ :

LEMMA 1. For every $n > 0$ we have :
$$\prod_{j=0}^{n-1}(X - Y^j) \equiv X^n - 1 \pmod{\Phi_n(Y)}. \qquad (2)$$

PROOF. Consider the polynomial
$$F(X,Y) = \prod_{j=0}^{n-1}(X - Y^j) - (X^n - 1).$$

We can write it in the form :
$$f_0(Y) + f_1(Y)X + f_2(Y)X^2 + \cdots + f_{n-1}(Y)X^{n-1},$$

with $f_k(Y) \in \mathbb{Z}[Y]$. When $\zeta$ is a primitive $n$-th root of unity, the numbers $\zeta^j$ run through all $n$-th roots of unity when $j$ takes on the values $0, 1, ..., n-1$, whence $F(X, \zeta)$ is the zero polynomial in $X$ for every primitive $n$-th root $\zeta$ of unity. This means that, for $k = 0$, 1, ..., $n - 1$, the polynomial $f_k(Y)$ is zero at every root of $\Phi_n(Y)$, which proves lemma 1. ∎

---

[*] Address : Leeuwerikenstraat 74, B-3001 Heverlee, Belgium.
 manuvti@hotmail.com.

Dividing both sides of formula (2) by $X - 1$ is possible and yields :

COROLLARY. For every $n > 1$ we have
$$\prod_{j=1}^{n-1}(X - Y^j) \equiv X^{n-1} + X^{n-2} + \cdots + X + 1 \ (\mathrm{mod}\ \Phi_n(Y)). \qquad (3)$$

PROOF OF THEOREM 1.
Let $m$ be a prime number. Then we have
$$\Phi_m(X) = \frac{X^m - 1}{X - 1}.$$
Putting $n = m$, $X = -1$, $Y = 2$ in (3), we obtain :
$$\prod_{j=1}^{m-1}(-1 - 2^j) \equiv 1\ (\mathrm{mod}\ 2^m - 1).$$
which implies (1) in case $m$ is an odd prime.

For the converse, suppose (1) holds for an $m \geq 2$. If $m$ would be even, say $m = 2w$, then the modulus of the congruence (1) would be divisible by $2^w + 1$, a divisor of the left hand of (1) but not of the right hand, which is impossible. Therefore, (1) implies that $m$ is odd. Assume $m$ composite. Let $c$ be a proper divisor of $m$ and $d = m/c$. Then we have, after multiplication of (1) by 2 (which extends the product in (1) to start with $j = 0$) and grouping the $m$ factors into $c$ blocks of $d$ factors, that
$$\prod_{i=0}^{m-1}(2^i + 1) \equiv \prod_{k=0}^{c-1}\left(\prod_{j=0}^{d-1}(2^{kd+j} + 1)\right) \equiv 2\ (\mathrm{mod}\ 2^m - 1).$$
From (*) we get, after replacing $n$ by $m$ and $X$ by 2, that $\Phi_d(2)$ divides $2^m - 1$, whence the former congruence holds also modulo $\Phi_d(2)$. By appropiate use of (*), we also see that $2^d - 1$ is divisible by $\Phi_d(2)$, whence $2^d \equiv 1\ (\mathrm{mod}\ \Phi_d(2))$. Thus
$$\prod_{j=0}^{d-1}(2^{kd+j} + 1) \equiv \prod_{j=0}^{d-1}(2^j + 1)\ (\mathrm{mod}\ \Phi_d(2)), \forall k = 0, 1, \ldots, c-1.$$
Now, setting $n = d$, $X = -1$ and $Y = 2$ in (2), we get (after changing signs)
$$\prod_{j=0}^{d-1}(2^j + 1) \equiv 2\ (\mathrm{mod}\ \Phi_d(2)).$$
From the last three congruences we obtain $2^c \equiv 2\ (\mathrm{mod}\ \Phi_d(2))$ or
$$2^{c-1} \equiv 1\ (\mathrm{mod}\ \Phi_d(2)). \qquad (**)$$

Recall that the multiplicative order of $a$ modulo $b$ is defined as the least positive integer $u$ such that $a^u \equiv 1\ (\mathrm{mod}\ b)$. If $(a,b) = 1$, then such $u$ allways exists and if $v$ has the property $a^v \equiv 1\ (\mathrm{mod}\ b)$ then $v$ is divisible by $u$. From the congruence (**) we can deduce that the multiplicative order $t$ of 2 modulo $\Phi_d(2)$ must be a divisor of $c - 1$. But, since $2^d - 1$ is divisible by $\Phi_d(2)$, $t$ should be a divisor of $d$. If we let $c$ be the least prime divisor of $m$, no prime divisor of $d$ can divide $c - 1$, making (**) impossible, which proves theorem 1. ∎

To obtain a generalization of our result, we follow a suggestion of Kilford in [1]. From (3), replacing $n$ by $m$, $X$ by $-1$ and $Y$ by an integer $M > 2$, we easily get :

$$\prod_{j=1}^{m-1}(M^j + 1) \equiv 1 \pmod{\frac{M^m - 1}{M - 1}}, \qquad (4)$$

when $m$ is an odd prime. Next we prove that (4) only holds when $m$ is an odd prime.

Suppose that (4) holds for some $m > 2$. If $m = 2w$, then $M^w + 1$ divides the modulus and the left hand of (4) which produces a contradiction. Thus, $m$ must be odd. Assume $m$ composite, say $m = c\,d$ with $c$ and $d > 1$. Then, we can deduce in the same way as in the proof of Theorem 1 that

$$\prod_{i=0}^{m-1}(M^i + 1) \equiv \prod_{k=0}^{c-1}\left(\prod_{j=0}^{d-1}(M^{kd+j} + 1)\right) \equiv 2 \pmod{\frac{M^m - 1}{M - 1}}.$$

Again by (*), this congruence also holds modulo $\Phi_d(M)$. In the same way as in the proof of Theorem 1 we can derive that

$$2^{c-1} \equiv 1 \pmod{\Phi_d(M)}. \qquad (***)$$

Unfortunately, we cannot continue using the argumentation at the end of Theorem 1. Nevertheless, we can rely on the following result :

LEMMA 2. For $M > 2$ and $d > 2$, we have $\Phi_d(M) > (M-1)^{\varphi(d)}$, where $\varphi(d)$ is Eulers totient function.

PROOF. When $d > 2$, the polynomial $\Phi_d(X)$ has no real roots, whence $\Phi_d(M)$ is positive. Thus :

$$\Phi_d(M) = |\Phi_d(M)| = \prod |M - \zeta|,$$

where the product is taken over all the primitive $d$-th roots of unity. The factor $|M - \zeta|$ is obviously strictly greater than $M - 1$, and there are $\varphi(d)$ such factors. This terminates the proof of lemma 2. ∎

It is now clear that the congruence (***) cannot hold when we may choose $d$ to be such that each prime divisor $p$ of $d$ is $\geq c$. Indeed, $\varphi(d)$ is then divisible by $p - 1$, which implies that $\Phi_d(M) > 2^{c-1}$ which prevents $2^{c-1} - 1$ to be divisible by $\Phi_d(M)$.

All this finally yields :

THEOREM 2.
Let $M$ be a natural number, stricly greater than 1. Then, the number $m > 2$ is prime if and only if

$$\prod_{j=1}^{m-1}(M^j + 1) \equiv 1 \pmod{\frac{M^m - 1}{M - 1}}.$$

ACKNOWLEDGEMENT

I wish to express my thanks to Prof. J. Denef for his help.